\documentclass[10pt]{article}
\usepackage{amsfonts}
\usepackage{amssymb,amsmath,amsthm, amsfonts}

\def\sqr#1#2{{\vcenter{\vbox{\hrule height.#2pt
              \hbox{\vrule width.#2pt height#1pt \kern#1pt \vrule width.#2pt}
          \hrule height.#2pt}}}}
\def\signed #1{{\unskip\nobreak\hfil\penalty50
          \hskip2em\hbox{}\nobreak\hfil#1
          \parfillskip=0pt \finalhyphendemerits=0 \par}}
\def\endpf{\signed {$\sqr69$}}
\def\sqr#1#2{{\vcenter{\vbox{\hrule height.#2pt
              \hbox{\vrule width.#2pt height#1pt \kern#1pt \vrule width.#2pt}
              \hrule height.#2pt}}}}
\def\signed #1{{\unskip\nobreak\hfil\penalty50
              \hskip2em\hbox{}\nobreak\hfil#1
              \parfillskip=0pt \finalhyphendemerits=0 \par}}
\def\endpf{\signed {$\sqr69$}}
\def\3n{\negthinspace \negthinspace \negthinspace }
\def\2n{\negthinspace \negthinspace }
\def\1n{\negthinspace }

\def\={\buildrel \triangle \over =}

\def\esssup{\mathop{\rm esssup}}

\def\sup{\mathop{\rm sup}}

\def\esssup{\hbox{\rm ess$\,$\rm sup$\,$}}

\def\|{\Big |}
\def\({\Big (}
\def\){\Big )}
\def\[{\Big[}
\def\]{\Big]}
\def\be{\begin{equation}}
\def\bel{\begin{equation}\label}
\def\ee{\end{equation}}
\def\bt{\begin{theorem}}
\def\bcd{\begin{condition}}
\def\ecd{\end{condition}}
\def\et{\end{theorem}}
\def\bc{\begin{corollary}}
\def\ec{\end{corollary}}
\def\bde{\begin{definition}}
\def\ede{\end{definition}}
\def\bl{\begin{lemma}}
\def\el{\end{lemma}}
\def\bp{\begin{proposition}}
\def\ep{\end{proposition}}
\def\bex{\begin{example}}
\def\eex{\end{example}}
\def\br{\begin{remark}}
\def\er{\end{remark}}
\def\ba{\begin{array}}
\def\ea{\end{array}}
\def\ed{\end{document}}

\def\square#1{\vbox{\hrule\hbox{\vrule height#1%
     \kern#1\vrule}\hrule}}
\def\rectangle#1#2{\vbox{\hrule\hbox{\vrule height#1%
     \kern#2\vrule}\hrule}}

\font\tenbb=msbm10 \font\sevenbb=msbm7 \font\fivebb=msbm5

\newfam\bbfam
\scriptscriptfont\bbfam=\fivebb \textfont\bbfam=\tenbb
\scriptfont\bbfam=\sevenbb

\newtheorem{lemma}{Lemma}[section]
\newtheorem{remark}{Remark}[section]
\newtheorem{example}{Example}[section]
\newtheorem{theorem}{Theorem}[section]
\newtheorem{corollary}{Corollary}[section]

\newtheorem{definition}{Definition}[section]
\newtheorem{proposition}{Proposition}[section]
\newtheorem{condition}{Condition}[section]

\makeatletter
   
   \@addtoreset{equation}{section}
\makeatother

\begin{document}

\title{Note on stochastic control problems related with general fully coupled forward-backward stochastic differential equations}
\author{Juan Li\\
{\small School of Mathematics and Statistics, Shandong University at Weihai, Weihai 264209, P. R. China.}\\
{\small{\it E-mail:  juanli@sdu.edu.cn}}}
\date{June 21, 2012}

\maketitle \noindent{\bf Abstract.}\hskip4mm  In this paper we study
stochastic optimal control problems of general fully coupled
forward-backward stochastic differential equations (FBSDEs). In Li and Wei~\cite{LW} the authors studied two cases of diffusion coefficients $\sigma$ of
FSDEs, in one case when $\sigma$\ depends on the control and does not depend on the second component of the solution $(Y, Z)$ of the BSDE, and in the other case $\sigma$ depends on $Z$ and doesn't depend on the control. Here we study the general case when $\sigma$ depends on both $Z$ and the control at the same time. The recursive cost functionals are defined by controlled general fully coupled FBSDEs, then the value functions are given by taking the essential supremum of the cost functionals over all admissible controls. We give the formulation of related generalized Hamilton-Jacobi-Bellman (HJB) equations, and prove the value function is its viscosity solution.

\medskip


\newpage

\section{\large{Introduction}}

\hskip1cm

Pardoux and Peng~\cite{PaPe1} first introduced nonlinear backward stochastic differential equations
(BSDEs) driven by a Brownian motion. Since then the theory of BSDEs develops very quickly, see
El Karoui, Peng and Quenez~\cite{ELPeQu}, Peng~\cite{Pe1},~\cite{Pe2},~\cite{Pe3}, etc. Associated with the BSDE theory, the theory of fully coupled
forward-backward stochastic differential equations (FBSDEs) develops also very quickly, refer to, Antonelli~\cite{An},
Cvitanic and Ma~\cite{CM}, Delarue~\cite{D}, Hu and Peng~\cite{HP}, Ma, Protter, and Yong~\cite{MPY}, Ma, Wu, Zhang, and Zhang~\cite{MWZZ}, Ma and Yong~\cite{MY}, Pardoux and
Tang~\cite{PaT}, Peng and Wu~\cite{PW}, Wu~\cite{W}, Yong~\cite{Y},~\cite{Y2}, and
Zhang~\cite{Z}, etc. For more details on fully coupled FBSDEs, the
reader is referred to the book of Ma and Yong~\cite{MY}; also refer to Li and Wei~\cite{LW}
and the references therein.

Pardoux and Tang~\cite{PaT} studied fully coupled FBSDEs (but without controls), and gave an existence result for viscosity
solutions of related quasi-linear parabolic PDEs, when the diffusion coefficient $\sigma$ of the
forward equation does not depend on the second component of the solution $(Y, Z)$ of the BSDE.  Wu and Yu~\cite{WY},~\cite{WY2} studied the case when the diffusion coefficient $\sigma$ of the forward equation depends on $Z$, but, the stochastic systems without controls. In Li and Wei~\cite{LW} they studied the optimal control problems of fully coupled FBSDEs. They studied two cases of diffusion coefficients $\sigma$ of FSDEs, that is, in one case $\sigma$\ depends on the control and does not depend on $Z$, and in the other case $\sigma$ depends on $Z$ and doesn't depend on the control. They use a new method to prove that the value functions are deterministic, satisfied the dynamic programming principle (DPP), and were viscosity solutions to the related generalized
Hamilton-Jacobi-Bellman (HJB) equations. The associated generalized HJB equations are related with algebraic equations when $\sigma$
depends on $Z$ and doesn't depend on the control. They generalized Peng's BSDE
method, and in particular, the notion of stochastic backward
semigroup in~\cite{Pe4}. When $\sigma$ depends on $Z$ makes the stochastic control much more
complicate and the related HJB equation is combined with an algebraic equation, which was inspired by Wu,
Yu~\cite{WY}. However, they use the continuation method combined with the fixed point theorem to prove
very clearly that the algebraic equation has a unique solution, and, moreover, also give the representation for this solution. On the other hand, they also prove some new basic estimates for fully coupled FBSDEs under the monotonic assumptions. In particular, they prove under the Lipschitz and linear growth conditions that fully coupled FBSDEs have a unique solution on the small time interval, if the Lipschitz constant of $\sigma$\ with respect to $z$ is sufficiently small. They also establish a generalized comparison theorem for such fully coupled FBSDEs. Here we want to study the general case, that is, when $\sigma$ depends on $Z$ and the control at the same time, for this general case what the associated HJB equations will be.  

Let us be more precise. We study a stochastic control problem related with fully coupled FBSDE. The cost functional is introduced by the
following fully coupled FBSDE: \be\label{ee1.1}
  \left \{
  \begin{array}{llll}
  dX_s^{t,x;u} & = & b(s,X_s^{t,x;u},Y_s^{t,x;u},Z_s^{t,x;u},u_s)ds +
          \sigma(s,X_s^{t,x;u},Y_s^{t,x;u},Z_s^{t,x;u},u_s) dB_s, \\
  dY_s^{t,x;u} & = & -f(s,X_s^{t,x;u},Y_s^{t,x;u},Z_s^{t,x;u},u_s)ds + Z_s^{t,x;u}dB_s, \ \ \ \ \ s\in[t,T], \\
   X_t^{t,x;u}& = & x,\\
   Y_T^{t,x;u} & = & \Phi(X_T^{t,x;u}),
   \end{array}
   \right.
  \ee
\\where $T>0$ is an arbitrarily fixed finite time horizon,
$B=(B_s)_{s\in[0,T]}$ is a d-dimensional standard Brownian motion,
and $u=(u_s)_{s\in[t,T]}$ is an admissible control. Precise
assumptions on the coefficients $b,\  \sigma, \ f,  \  \Phi$ are given in the next section. Under
 our assumptions, (\ref{ee1.1}) has a unique solution
  $(X_s^{t,x;u},Y_s^{t,x;u},Z_s^{t,x;u})_{s\in[t,T]}$ and the cost
  functional is defined by \be\label{ee1.2} J(t,x;u)=Y_t^{t,x;u}. \ee

\noindent We define the value function of our stochastic control
problems as follows: \be\label{ee1.3} W(t,x):=\esssup_{u\in\mathcal
{U}_{t,T}}J(t,x;u). \ee The objective of our paper is to investigate
this value function. The main results of the paper state that $W$
is deterministic (Proposition 2.1), continuous viscosity solution
of the associated HJB equations (Theorem 3.1). The associated HJB
equation is combined with an algebraic equation as follows:
\be\label{ee1.5}
 \left \{\begin{array}{ll}
 &\!\!\!\!\! \frac{\partial }{\partial t} W(t,x) +  H_V(t, x, W(t,x))=0,
   \\
&\!\!\!\!\!V(t,x,u)=DW(t,x).\sigma(t,x,W(t,x),V(t,x,u),u),\hskip 0.5cm   (t,x)\in [0,T)\times {\mathbb{R}}^n , u\in U,\\
 &\!\!\!\!\!  W(T,x) =\Phi (x),\hskip 0.5cm   x\in {\mathbb{R}}^n.
 \end{array}\right.
\ee
 In this case
$$\begin{array}{lll}
H_V(t, x, W(t,x))&=& \sup\limits_{u \in U}\{DW.b(t, x, W(t,x),
V(t,x,u), u)+\frac{1}{2}tr(\sigma\sigma^{T}(t,
x, W(t,x), V(t,x,u),u)D^2W(t,x))\\
 &&
  +f(t, x, W(t,x), V(t,x,u), u)\},
 \end{array}$$
where $t\in [0, T], x\in{\mathbb{R}}^n.$

Our paper is organized as follows: Section 2 introduces the framework of the stochastic
control problems. In Section 3, we prove that $W$ is a viscosity solution of the
associated HJB equation described above.

\section{ {\large Framework}}
   \hskip1cm
Let $(\Omega, {\cal{F}},$ P$)$ be the classical Wiener space, where $\Omega$
is the set of continuous functions from [0, T] to ${\mathbb{R}}^d$
starting from 0 ($\Omega= C_0([0, T];{\mathbb{R}}^d)$), $ {\cal{F}}
$ is the completed Borel $\sigma$-algebra over $\Omega$, and P is the
Wiener measure. Let B be the canonical process:
$B_s(\omega)=\omega_s,\ s\in [0, T],\ \omega\in \Omega$. We denote
by ${\mathbb{F}}=\{{\mathcal{F}}_s,\ 0\leq s \leq T\}$\ the natural
filtration generated by $\{B_t\}_{t\geq0}$ and augmented by all
$P\mbox{-}$null sets, i.e., $$ {\cal{F}}_s=\sigma\{B_r,r\leq s\}\vee
\mathcal {N}_P,\ \ \ \ s\in [0,T],$$ where $\mathcal {N}_P$ is the
set of all $P\mbox{-}$null subsets and $T$ is a fixed real time
horizon.
 For any
    $n\geq 1,$\ $|z|$ denotes the Euclidean norm of $z\in
    {\mathbb{R}}^{n}$. We introduce the following two spaces of processes which will be used later:

 ${\cal{S}}^2(t_0, T; {\mathbb{R}^n})$\ is the set of $\mathbb{R}^n\mbox{-valued}\ \mathbb{F}\mbox{-adapted
    continuous process}\ (\psi_t)_{t_0\leq t\leq T} $ which also satisfies $E[\sup\limits_{t_0\leq t\leq T}| \psi_{t} |^2]< +\infty; $ \    ${\cal{H}}^{2}(t_0,T;{\mathbb{R}}^{n})$\ is the set of ${\mathbb{R}}^{n}\mbox{-valued}\ \mathbb{F}
    \mbox{-progressively measurable process}\ (\psi_t)_{t_0\leq t\leq T}$ which also satisfies $\parallel\psi\parallel^2=E[\int^T_{t_0}| \psi_t| ^2dt]<+\infty; $\ where $t_0\in [0,T].$

First we introduce the setting for stochastic optimal control problems. We suppose that
the control state space U is a compact metric space. ${\mathcal{U}}$
is the set of all U-valued $ {{\mathbb {F}} }$-progressively
measurable processes. If $u\in \mathcal {U}$, we call $u$ an
admissible control.

For a given admissible control $u(\cdot)\in {\mathcal{U}}$, we
regard $t$ as the initial time and $\zeta \in L^2 (\Omega
,{\mathcal{F}}_t, P;{\mathbb{R}}^n)$ as the initial state. We
consider the following fully coupled forward-backward stochastic
control system

\be\label{3.1}
  \left \{
  \begin{array}{llll}
  dX_s^{t,\zeta,u} & = & b(s,X_s^{t,\zeta;u},Y_s^{t,\zeta;u},Z_s^{t,\zeta;u},u_s)ds +
          \sigma(s,X_s^{t,\zeta;u},Y_s^{t,\zeta;u},Z_s^{t,\zeta;u},u_s) dB_s, \\
  dY_s^{t,\zeta;u} & = & -f(s,X_s^{t,\zeta;u},Y_s^{t,\zeta;u},Z_s^{t,\zeta;u},u_s)ds + Z_s^{t,\zeta;u}dB_s, \ \ \ \ \ \ \ \ \ \  s\in[t,T],\\
   X_t^{t,\zeta;u}& = & \zeta,\\
   Y_T^{t,\zeta;u} & = & \Phi(X_T^{t,\zeta;u}),
   \end{array}
   \right.
  \ee
  where the deterministic mappings

  $b: [0,T] \times \mathbb{R}^n \times  \mathbb{R} \times  \mathbb{R}^d  \times U
  \longrightarrow \mathbb{R}^n, \ \ \ \ \ \ \ \ \ \ \ \  \  \sigma:  [0,T] \times \mathbb{R}^n \times  \mathbb{R} \times  \mathbb{R}^d \times U
  \longrightarrow \mathbb{R}^{n\times d}, $

  $f:  [0,T] \times \mathbb{R}^n \times  \mathbb{R} \times  \mathbb{R}^d \times U
  \longrightarrow \mathbb{R},\ \  \ \ \ \ \ \ \ \ \ \ \ \ \Phi: \mathbb{R}^n
  \longrightarrow \mathbb{R}$\\
  are continuous in $(t, u)\in [0, T]\times U$. In this paper
we use the usual inner product and the Euclidean norm in
$\mathbb{R}^n, \mathbb{R}^m$ and $\mathbb{R}^{m\times d},$
respectively. Given an $m \times n$ full-rank matrix G, we define:
$$\lambda= \
            \left(\begin{array}{c}
            x\\
            y\\
            z
            \end{array}\right)
            \ , \ \ \ \ \ \ \ \ \ \
A(t,\lambda)= \
            \left(\begin{array}{c}
            -G^Tf\\
            Gb\\
            G\sigma
            \end{array}\right)(t,\lambda),$$where $G^T$ is the transposed matrix of $G$.

We assume that

(B1) (i) $A(t,\lambda)$ is uniformly Lipschitz  with respect to
$\lambda$, and for any $\lambda$, $A(\cdot,\lambda)\in $

\hskip1.4cm${\cal{H}}^{2}(0,T;{\mathbb{R}}^{n}\times{\mathbb{R}}^{m}\times{\mathbb{R}}^{m\times
d});$

 \ \ \ \ \ \ \  (ii) $\Phi(x)$ is uniformly Lipschitz
with respect to $x\in \mathbb{R}^n$, and for any $x\in
\mathbb{R}^n,\Phi(x)\in$

 \hskip1.4cm$ L^2(\Omega,\mathcal
{F}_T,\mathbb{R}^m).$

\noindent The following monotonicity conditions are also necessary:

 (B2) (i) $\langle A(t,\lambda)-A(t,\overline{\lambda}),\lambda-\overline{\lambda} \rangle \leq -\beta_1|G\widehat{x}|^2-\beta_2(|G^T \widehat{y}|^2+|G^T
 \widehat{z}|^2),$

\hskip0.8cm (ii) $ \langle
\Phi(x)-\Phi(\overline{x}),G(x-\overline{x}) \rangle \geq
\mu_1|G\widehat{x}|^2,\ \ \widehat{x}=x-\bar{x},\
\widehat{y}=y-\bar{y},\ \widehat{z}=z-\bar{z}$,
\\where $\beta_1,\ \beta_2,\ \mu_1$ are nonnegative constants with
$\beta_1 + \beta_2>0,\ \beta_2 + \mu_1>0$. Moreover, we have
$\beta_1>0,\ \mu_1>0 \ (\mbox{resp., }\beta_2>0)$, when $m>n$ (resp.,
$m<n$).

\br \rm{(B3)-(ii)} $ \langle
\Phi(x)-\Phi(\overline{x}),G(x-\overline{x}) \rangle \geq
0.$
\er

The coefficients satisfy the assumptions (B1) and (B2), and
  also\\
(B4) there exists a constant $K\geq 0$ such that, for all $t\in [0, T],\ u\in U,\ x_1,\
x_2\in\mathbb{R}^n, \ y_1,\ y_2\in\mathbb{R},\ z_1,\
z_2\in\mathbb{R}^d,$
$$|l(t,x_1,y_1,z_1,u)-l(t,x_2,y_2,z_2,u)|\leq K(|x_1-x_2|+|y_1-y_2|+|z_1-z_2|),$$
\ \ \ \ \ \ $l=b, \sigma, f$, respectively, and $|\Phi(x_1)-\Phi(x_2)|\leq K|x_1-x_2|$.
\br Under our assumptions, it is obvious that there exists a constant $C\geq 0$ such that,
$$|b(t,x,y,z,u)|+|\sigma(t,x,y,z,u)|+|f(t,x,y,z,u)|+|\Phi(x)|\leq C(1+|x|+|y|+|z|),$$
for all $(t,x,y,z,u)\in [0,T]\times
\mathbb{R}^n\times\mathbb{R}\times\mathbb{R}^d\times U.$ \er

  Hence, for any $u(\cdot) \in \mathcal
  {U},$ from Lemma 2.4 in~\cite{LW}, FBSDE (\ref{3.1}) has a unique solution.

   From Proposition 6.1 in~\cite{LW}, there exists $C \in \mathbb{R}^+$ such that, for any $t \in
  [0,T]$, $\zeta, \zeta' \in L^2(\Omega,\mathcal
  {F}_t,P;\mathbb{R}^n),$ $u(\cdot) \in \mathcal
  {U},$ we have, \mbox{P-a.s.}: \be\label{3.2}
  \begin{array}{llll}
 \mbox{(i)}&E[\sup\limits_{t\leq s\leq T}|X_s^{t,\zeta;u}-X_s^{t,\zeta';u}|^2+\sup\limits_{t\leq s\leq T}|Y_s^{t,\zeta;u}-Y_s^{t,\zeta';u}|^2+\int_t^T|Z_s^{t,\zeta;u}-Z_s^{t,\zeta';u}|^2ds\mid\mathcal {F}_t]  \leq  C|\zeta - \zeta'|^2, \\
 \mbox{(ii)}& E[\sup\limits_{t\leq s\leq T}|X_s^{t,\zeta;u}|^2+\sup\limits_{t\leq s\leq T}|Y_s^{t,\zeta;u}|^2+ \int_t^T|Z_s^{t,\zeta;u}|^2ds \mid\mathcal {F}_t] \leq   C(1 +|\zeta|^2). \\
  \end{array}
  \ee
  Therefore, we get
 \be\label{3.3} \begin{array}{llll}&&{\rm(i)}\ \ \
  |Y_t^{t,\zeta;u}| \leq C(1+|\zeta|), \mbox{P-a.s.};\\
  &&{\rm(ii)}\ \ |Y_t^{t,\zeta;u}-Y_t^{t,\zeta';u}| \leq C|\zeta - \zeta'|, \mbox{P-a.s.}  \end{array}\ee

We now introduce the subspaces of admissible controls. An admissible
control process $u=(u_r)_{r\in [t,s]}$ on $[t,s]$ is an
$\mathbb{F}$-progressively measurable, $U$-valued process. The set
of all admissible controls on $[t,s]$ is denote $\mathcal
{U}_{t,s},\ t\leq s\leq T.$

For a given process $u(\cdot) \in \mathcal {U}_{t,T}$, we define the
associated cost functional as follows:
  \be\label{3.4} J(t,x;u):=Y_s^{t,x;u}\mid_{s=t}, \ \ (t,x)\in[0,T] \times
  \mathbb{R}^n,\ee where the process $Y^{t,x;u}$ is defined by FBSDE
  (\ref{3.1}).

From Theorem 6.1 in~\cite{LW} we have, for any $t \in
  [0,T]$ and $\zeta \in L^2(\Omega,\mathcal {F}_t,P;\mathbb{R}^n),$
 \be\label{3.5} J(t,\zeta;u)=Y_t^{t,\zeta;u},  \mbox{ P-a.s.} \ee
For $ \zeta = x \in \mathbb{R}^n,$ we define the value
  function as \be\label{3.6} W(t,x) := \esssup_{u\in\mathcal {U}_{t,T}}J(t,x;u). \ee
  \br From the assumptions (B1) and (B2), the value
  function $W(t,x)$ as the essential supremum over a family of ${\cal F}_t$-measurable random variables is well defined and it is a bounded ${\mathcal {F}}_t$-measurable random variable too. But it turns out to be even deterministic. \er
  In fact, inspired by the method in Buckdahn and Li~\cite{BL}, we can prove
  that $W$ is deterministic.
  \bp Under the assumptions (B1) and (B2), for any $(t,x) \in [0,T] \times \mathbb{R}^n,$ $W(t,x)$ is a deterministic function in the sense that $W(t,x) = E[W(t,x)],
  \mbox{P-a.s.}$
  \ep
The proof can be consulted in Li and Wei~\cite{LW}, Proposition 3.1.\\ \endpf

 From (\ref{3.3}) and (\ref{3.6}) we get the following property of the value function
$W(t,x)$:

\bl There exists a constant $C>0$ such that, for all $0 \leq t \leq
T,\ x,x' \in \mathbb{R}^n,$ \be \begin{array}{llll}&& {\rm(i)} \ \
|W(t,x)-W(t,x')| \leq C|x-x'|; \\
&&  {\rm(ii)}\ \ |W(t,x)| \leq C(1+|x|). \end{array}\ee \el

\bl Under the assumptions (B1) and (B2), the cost functional
$J(t,x;u),$\ for any $u\in \mathcal {U}_{t,T}$, and the value
function $W(t,x)$ are monotonic in the following sense: for each
$x,\ \bar{x} \in \mathbb{R}^n,$ $t\in [0,T],$
$$\begin{array}{llll}
&&{\rm(i)}\ \  \langle
J(t,x;u)-J(t,\bar{x};u),\ G(x-\bar{x})\rangle\geq 0,\ \mbox{P-a.s.};\\
&&{\rm(ii)}\ \ \langle
W(t,x)-W(t,\bar{x}),\ G(x-\bar{x})\rangle \geq 0.\end{array}$$  \el
 For the proof the reader is referred to Lemma 3.3 in Li and Wei~\cite{LW}. 
\br (1) From (B2)-(i) we see that if $\sigma$\ doesn't depend on $z$, then $\beta_2=0$. Furthermore, we assume that:\\ $\begin{array}{llll}
&&{\rm{(B5)}}\  \mbox{the Lipschitz constant}\ L_\sigma\geq 0 \ \mbox{of}\ \sigma\   \mbox{with respect to}\ z\ \mbox{is sufficiently small, i.e., there exists some}\\
&&\ \ \ \ \ \ L_\sigma\geq 0\ \mbox{small enough such that}, \ \mbox{for all}\ t\in[0, T],\ u\in U,\ x_1,\ x_2\in\mathbb{R}^n,\ y_1,\ y_2\in\mathbb{R},\ z_1,\ z_2\in\mathbb{R}^d,\\
&&\ \ \ \ \ \ |\sigma(t,x_1,y_1,z_1,u)-\sigma(t,x_2,y_2,z_2,u)|\leq K(|x_1- x_2|+|y_1-y_2|)+L_\sigma|z_1-z_2|.
\end{array}$\\
(2) On the other hand, notice that when $\sigma$\ doesn't depend on $z$ it's obvious that (B5) always holds true.
\er

 The notation of stochastic backward
  semigroup was first introduced by Peng~\cite{Pe4} and was applied to
  prove the DPP for stochastic control problems. Now we discuss a generalized DPP for our stochastic optimal control problem (\ref{3.1}), (\ref{3.6}). For this we have to adopt Peng's notion of stochastic backward
  semigroup, and to define the family of (backward) semigroups associated with FBSDE (\ref{3.1}).

    For given initial data $(t,x)$, a real $0<\delta \leq
    T-t,$ an admissible control process $u(\cdot)\ \in \ \mathcal
    {U}_{t,t+\delta}$ and a real-valued ${\cal F}_{t+\delta}\otimes {\cal B}(\mathbb{R}^n)$-measurable random function $\Psi: \Omega\times \mathbb{R}^n\rightarrow \mathbb{R}$, such that (B2)-(ii) holds, we put
    $$G_{s,t+\delta}^{t,x;u}[\Psi(t+\delta, \widetilde{X}_{t+\delta}^{t,x;u})]:=\widetilde{Y}_s^{t,x;u}, \  s \in
    [t,t+\delta],$$
    where $(\widetilde{X}_s^{t,x;u},\widetilde{Y}_s^{t,x;u},\widetilde{Z}_s^{t,x;u})_{t \leq s \leq
    t+\delta}$ is the solution of the following FBSDE with time
    horizon $t+\delta$:
\be\label{3.8}
  \left \{
  \begin{array}{llll}
d\widetilde{X}_s^{t,x;u} & = & b(s,\widetilde{X}_s^{t,x;u},\widetilde{Y}_s^{t,x;u},\widetilde{Z}_s^{t,x;u},u_s)ds + \sigma(s,\widetilde{X}_s^{t,x;u},\widetilde{Y}_s^{t,x;u},\widetilde{Z}_s^{t,x;u},u_s)dB_s, \\
 d\widetilde{Y}_s^{t,x;u} & = & -f(s,\widetilde{X}_s^{t,x;u},\widetilde{Y}_s^{t,x;u},\widetilde{Z}_s^{t,x;u},u_s)ds + \widetilde{Z}_s^{t,x;u}dB_s, \ \ \ \ \ s\in [t,t+\delta],\\
  \widetilde{X}_t^{t,x;u} & = & x,\\
  \widetilde{Y}_{t+\delta}^{t,x;u} & = & \Psi(t+\delta, \widetilde{X}_{t+\delta}^{t,x;u}).
   \end{array}
   \right.
  \ee
\br {\rm(i)} From Lemmas 2.4 and 2.5 in~\cite{LW} we know that, if $\Psi$\ doesn't depend on $x$, FBSDE (\ref{3.8}) has a unique solution
$(\widetilde{X}^{t,x;u},\widetilde{Y}^{t,x;u},\widetilde{Z}^{t,x;u}).$\\
{\rm(ii)} We also point out that if $\Psi$\ is Lipschitz with respect to $x$, FBSDE (\ref{3.8}) can be also solved under the assumptions (B4) and (B5) on the small interval $[t, t+\delta]$, for any $0\leq \delta\leq \delta_0,$\ where the small parameter $\delta_0>0$\ is independent of $(t, x)$\ and the control $u$; see Proposition 6.4 in~\cite{LW}. \er
Since $\Phi$\ satisfies (B2)-(ii), the solution $(X^{t,x;u},Y^{t,x;u},Z^{t,x;u})$ of FBSDE
(\ref{3.1}) exists and we get
$$ G_{t,T}^{t,x;u}[\Phi(X_T^{t,x;u})] =
G_{t,t+\delta}^{t,x;u}[Y_{t+\delta}^{t,x;u}]. $$ Moreover, we have
\be\label{3.9} J(t,x;u)=Y_t^{t,x;u}=G_{t,T}^{t,x;u}[\Phi(X_T^{t,x;u})]
=G_{t,t+\delta}^{t,x;u}[Y_{t+\delta}^{t,x;u}]=G_{t,{t+\delta}}^{t,x;u}[J(t+\delta,X_{t+\delta}^{t,x;u};u)].\ee
 \bt Under the assumptions (B1), (B2), (B4) and (B5), the value function $W(t,x)$
 satisfies the following DPP: there exists a sufficiently small $\delta_0>0$, such that for any $0\leq \delta \leq \delta_0,\ t\in [0, T-\delta],\ x \in \mathbb{R}^n,$
$$W(t,x)=\esssup_{u\in \mathcal
{U}_{t,t+\delta}}G_{t,{t+\delta}}^{t,x;u}[W(t+\delta,\widetilde{X}_{t+\delta}^{t,x;u})].$$
\et
The proof refers to Theorem 3.1 in Li and
Wei~\cite{LW}.

 Notice that from the definition of our stochastic backward
  semigroup we know that $$G_{s,t+\delta}^{t,x;u}[W(t+\delta,\widetilde{X}_{t+\delta}^{t,x;u})]=\widetilde{Y}_s^{t,x;u}, \  s \in
    [t,t+\delta],\ u(\cdot)\ \in \ \mathcal
    {U}_{t,t+\delta},$$
    where $(\widetilde{X}_s^{t,x;u},\widetilde{Y}_s^{t,x;u},\widetilde{Z}_s^{t,x;u})_{t \leq s \leq
    t+\delta}$ is the solution of the following FBSDE with time horizon $t+\delta$:
\be\label{3.10}
  \left \{
  \begin{array}{llll}
d\widetilde{X}_s^{t,x;u} & = & b(s,\widetilde{X}_s^{t,x;u},\widetilde{Y}_s^{t,x;u},\widetilde{Z}_s^{t,x;u},u_s)ds + \sigma(s,\widetilde{X}_s^{t,x;u},\widetilde{Y}_s^{t,x;u},\widetilde{Z}_s^{t,x;u},u_s)dB_s, \\
 d\widetilde{Y}_s^{t,x;u} & = & -f(s,\widetilde{X}_s^{t,x;u},\widetilde{Y}_s^{t,x;u},\widetilde{Z}_s^{t,x;u},u_s)ds + \widetilde{Z}_s^{t,x;u}dB_s, \ \ \ \ \ s\in [t,t+\delta],\\
  \widetilde{X}_t^{t,x;u} & = & x,\\
  \widetilde{Y}_{t+\delta}^{t,x;u} & = & W(t+\delta,\widetilde{X}_{t+\delta}^{t,x;u}).
   \end{array}
   \right.
  \ee
From Proposition 6.4 in~\cite{LW} there exists a sufficiently small $\delta_0>0$, such that for any $0\leq \delta \leq \delta_0,$\ the above equation (\ref{3.10}) has a unique solution $(\widetilde{X}^{t,x;u},\widetilde{Y}^{t,x;u},\widetilde{Z}^{t,x;u})$\ on the time interval $[t, t+\delta]$.

From Lemma 2.2, we get the Lipschitz property of the value function $W(t,x)$ in $x$, uniformly in $t$. Now from Theorem 2.1 we can get the continuity property of $W(t,x)$ in $t$ and to conclude:

 \bt  Under the assumptions (B1), (B2), (B4) and (B5) $W(t,x)$ is continuous in $t$.\et
  Its proof can be found in Li and
Wei~\cite{LW}, Theorem 3.2. 
 \section{\large Viscosity Solutions of HJB Equations}
  \hskip1cm
In this section we show that the value function $W(t, x)$ defined in
(\ref{3.6}) is a viscosity solution of the corresponding HJB equation. For
this we use Peng's BSDE
 approach~\cite{Pe4} developed for stochastic control problems of decoupled
 FBSDEs.

 Let us consider equation (\ref{3.1}): \be\label{4.1}
  \left \{
  \begin{array}{llll}
  dX_s^{t,x;u} & = & b(s,X_s^{t,x;u},Y_s^{t,x;u},Z_s^{t,x;u},u_s)ds +
          \sigma(s,X_s^{t,x;u},Y_s^{t,x;u},Z_s^{t,x;u},u_s) dB_s, \\
  dY_s^{t,x;u} & = & -f(s,X_s^{t,x;u},Y_s^{t,x;u},Z_s^{t,x;u},u_s)ds + Z_s^{t,x;u}dB_s, \ \ \ \  s\in [t,T],\\
   X_t^{t,x;u}& = & x,\\
   Y_T^{t,x;u} & = & \Phi(X_T^{t,x;u}).
   \end{array}
   \right.
  \ee
The related HJB equation is the following PDE combined with the
algebraic equation: \be\label{4.2}
 \left \{\begin{array}{ll}
 &\!\!\!\!\! \frac{\partial }{\partial t} W(t,x) +  H_V(t, x, W(t,x))=0,
   \\
&\!\!\!\!\!V(t,x,u)=DW(t,x).\sigma(t,x,W(t,x),V(t,x,u),u),\hskip 0.5cm   (t,x)\in [0,T)\times {\mathbb{R}}^n ,\ u\in U,\\
 &\!\!\!\!\!  W(T,x) =\Phi (x),\ \ \ \ x\in \mathbb{R}^n,
 \end{array}\right.
\ee
 where
$$\begin{array}{lll}
 H_V(t, x, W(t,x))= \sup\limits_{u \in U}\{ DW(t,x).b(t, x,
W(t,x), V(t,x,u), u)+f(t, x, W(t,x), V(t,x,u), u)\\
 \hskip 3.5cm +\frac{1}{2}tr(\sigma\sigma^{T}(t,
x, W(t,x), V(t,x,u),u)D^2W(t,x))\},\ t\in [0, T],\ x\in{\mathbb{R}}^n.
 \end{array}$$

We give the definition of viscosity solution for this kind of
PDE. For more details on viscosity solution refer to Crandall, Ishii, Lions~\cite{CIL}.

\bde\mbox{ } A real-valued
continuous function $W\in C([0,T]\times {\mathbb{R}}^n )$ is called \\
  {\rm(i)} a viscosity subsolution of equation (\ref{4.2}) if $W(T,x) \leq \Phi (x),\ \mbox{for all}\ x \in
  {\mathbb{R}}^n$, and if for all functions $\varphi \in C^3_{l, b}([0,T]\times
  {\mathbb{R}}^n)$ and for all $(t,x) \in [0,T) \times {\mathbb{R}}^n$ such that $W-\varphi $\ attains
  a local maximum at $(t, x)$,
 $$\left \{\begin{array}{ll}
 &\!\!\!\!\! \frac{\partial \varphi}{\partial t} (t,x) + H_{\psi}(t,x,\varphi(t,x)) \geq
 0,\\
&\!\!\!\!\!\mbox{where }\psi \mbox{ is the unique solution of the
following algebraic
equation:}\\
&\!\!\!\!\!\psi(t,x,u)=D\varphi(t,x).\sigma(t,x,\varphi(t,x),\psi(t,x,u),u),\ u\in U.
\end{array}\right.
$$
 \noindent{\rm(ii)} a viscosity supersolution of equation (\ref{4.2})
if $W(T,x) \geq \Phi (x), \mbox{for all}\ x \in
  {\mathbb{R}}^n$, and if for all functions $\varphi \in C^3_{l, b}([0,T]\times
  {\mathbb{R}}^n)$ and for all $(t,x) \in [0,T) \times {\mathbb{R}}^n$\ such that $W-\varphi $\ attains
  a local minimum at $(t, x)$,
$$\left \{\begin{array}{ll}
 &\!\!\!\!\! \frac{\partial \varphi}{\partial t} (t,x) + H_{\psi}(t,x,\varphi(t,x)) \leq
 0,\\
 &\!\!\!\!\!\mbox{where }\psi \mbox{ is the unique solution of the following algebraic
equation:}\\
&\!\!\!\!\!
\psi(t,x,u)=D\varphi(t,x).\sigma(t,x,\varphi(t,x,u),\psi(t,x),u),\ u\in U.
\end{array}\right.
$$
 \noindent{\rm(iii)} a viscosity solution of equation (\ref{4.2}) if it is both a viscosity sub- and supersolution of equation (\ref{4.2}). \ede
\br When $\sigma$\ depends on $z$ we need the test function $\varphi$ in Definition 3.1 satisfies the monotonicity condition {\rm(B3)-(ii)} and also the following technical assumptions:
\\
(B6) $\beta_2>0$;\\
(B7) $G\sigma(s, x, y, z)$ is continuous in $(s, u)$, uniformly with respect to $(x, y, z)\in {\mathbb{R}}\times {\mathbb{R}}^n\times {\mathbb{R}}^d$.\er
 \bt  Under
the assumptions (B1), (B2), (B4), (B5), (B6) and (B7), the value function $W$
is a viscosity solution of (\ref{4.2}).\et

We have the following important Representation Theorem for the
solution of the algebraic equation.

\bp For any $s\in [0, T],\ \zeta\in\mathbb{R}^d,\ y\in \mathbb{R},\
\bar{x}\in \mathbb{R}^n,\ u\in U$, there exists a unique $z$ such that
$z=\zeta+D\varphi(s,\bar{x})\sigma(s,\bar{x},y+\varphi(s,\bar{x}),z,u).
$ That means, the solution $z$ can be written as $z=h(s, \bar{x},y,
\zeta,u)$, where the function $h$\ is Lipschitz with respect to $\ y,
\ \zeta,$ and $|h(s, \bar{x},y, \zeta,u)|\leq
C(1+|\bar{x}|+|y|+|\zeta|).$ The constant $C$\ is independent of
$s,\ \bar{x},\ y, \ \zeta,\ u$. And $z=h(s, \bar{x},y, \zeta,u)$ is
continuous with respect to $(s,u)$. \ep

 \noindent \textbf{Proof}. It's obvious that we only need to consider the case when $\sigma$\ depends on $z$. Then we have (B6) and (B7) hold. Furthermore, following the proof of Proposition 4.1 in~\cite{LW} we can prove it.\endpf

 \noindent \textbf{Proof of Theorem 3.1}. Following the proof of Theorem 4.2 in~\cite{LW} we can prove it.\endpf

\section*{Acknowledgments}  Juan Li gratefully acknowledges financial support by the NSF of P.R.China (Nos. 10701050, 11071144), Shandong
Province (Nos. BS2011SF010), SRF for ROCS (SEM).

\end{document}